\documentclass{ifacconf}

\usepackage{graphicx}      % include this line if your document contains figures
\usepackage{natbib}        % required for bibliography

\usepackage{amsfonts}
\usepackage{amsmath}
\usepackage{amssymb}
\usepackage[mathscr]{eucal}
\usepackage{amsfonts}
\usepackage{amsmath}
\usepackage{amssymb}
\def\begcen{\begin{center}}
	\def\endcen{\end{center}}

\newcommand{\col}{ \mbox{col} }

\def\caly{{\cal Y}}
\def\calg{{\cal G}}

\def\call{{\cal L}}
\def\calw{{\cal W}}

\def\hal{{1 \over 2}}

\def\liminf{\lim_{t \to \infty}}

\def\L2{{\cal L}_2}
\def\L2e{{\cal L}_{2e}}

\def\rea{\mathbb{R}}

\def\sign{\mbox{sign}}
\def\adj{\mbox{adj}}

\def\x{{x}}

\def\begequarr{\begin{eqnarray}}
	\def\endequarr{\end{eqnarray}}
\def\begequarrs{\begin{eqnarray*}}
	\def\endequarrs{\end{eqnarray*}}
\def\begarr{\begin{array}}
	\def\endarr{\end{array}}
\def\begequ{\begin{equation}}
	\def\endequ{\end{equation}}
\def\lab{\label}
\def\begdes{\begin{description}}
	\def\enddes{\end{description}}
\def\begenu{\begin{enumerate}}
	\def\begite{\begin{itemize}}
		\def\endite{\end{itemize}}
	\def\endenu{\end{enumerate}}

\def\lef[{\left[\begin{array}}
	\def\rig]{\end{array}\right]}

\def\begcen{\begin{center}}
	\def\endcen{\end{center}}
\def\begrem{\begin{remark}\rm}
	\def\endrem{\end{remark}}
\def\begassum{\begin{assumption}}
	\def\endassum{\end{assumption}}
\def\begassums{\begin{assumption*}}
	\def\endassums{\end{assumption*}}
\def\begassu{\begin{ass}}
	\def\endassu{\end{ass}}
\def\beglem{\begin{lemma}}
	\def\endlem{\end{lemma}}
\def\begcor{\begin{corollary}}
	\def\endcor{\end{corollary}}
\def\begfac{\begin{fact}}
	\def\endfac{\end{fact}}

\def\liminf{\lim_{t \to \infty}}

\def\L2e{{\cal L}_{2e}}

\def\rea{\mathbb{R}}
\def\intnum{\mathbb{Z}}
\def\sign{\mbox{sign}}

\def\adj{\mbox{adj}}
\def\col{\mbox{col}}
\def\hal{{1 \over 2}}

\def\intnum{\mathbb{Z}}

\def\begsubequ{\begin{subequations}}
	\def\endsubequ{\end{subequations}}
\newcommand\dif{\mathrel{\overset{\makebox[0pt]{\mbox{\normalfont\tiny ${d \over dt}$}}}{\Longrightarrow}}}

\newcommand\hhh{\mathrel{\overset{\makebox[0pt]{\mbox{\normalfont\tiny $\times h_3$}}}{\Longrightarrow}}}
\def\begmat#1{\begin{bmatrix}#1\end{bmatrix}}
\def\begali#1{\begin{align}{#1}\end{align}}
\def\begalis#1{\begin{align*}{#1}\end{align*}}
\newtheorem{lemma}{Lemma}
\newtheorem{proposition}{Proposition}

\usepackage{xcolor}

\begin{document}
\begin{frontmatter}

\title{Parameter Estimation of Two Classes of Nonlinear Systems with Non-separable Nonlinear Parameterizations\thanksref{footnoteinfo}} 

\thanks[footnoteinfo]{This paper is supported by the Ministry of Science and Higher Education of Russian Federation, passport of goszadanie no. 2019-0898. This work is part of the project MAFALDA (PID2021-126001OB-C31) and MASHED (TED2021-129927B-I00) funded by MCIN/ AEI/10.13039/501100011033 and by the European Union Next GenerationEU/PRTR}

\author[First]{Romeo Ortega} 
\author[Second]{Alexey Bobtsov} 
\author[Third]{Ramon Costa-Castell\'o}
\author[Second]{Nikolay Nikolaev}

\address[First]{Departamento Acad\'{e}mico de Sistemas Digitales, ITAM, Ciudad de M\'exico, M\'{e}xico (e-mail: romeo.ortega@itam.mx)}
\address[Second]{Department of Control Systems and Robotics, ITMO University, Saint-Petersburg, Russia (e-mail: bobtsov@mail.ru, nikona@yandex.ru)}
\address[Third]{Universitat Politècnica de Catalunya (UPC), Spain (e-mail: ramon.costa@upc.edu)}

\begin{abstract}                % Abstract of not more than 250 words.
In this paper we address the challenging problem of designing globally convergent estimators for the parameters of nonlinear systems containing a {\em non-separable} exponential nonlinearity. This class of terms appears in many practical applications, and none of the existing parameter estimators is able to deal with them in an efficient way. The proposed estimation procedure is illustrated with two modern applications: fuel cells and human musculoskeletal dynamics. The procedure does not assume that the parameters live in known compact sets, that the nonlinearities satisfy some Lipschitzian properties, nor rely on injection of high-gain or the use of complex, computationally demanding methodologies. Instead, we propose to design a classical on-line estimator whose dynamics is described by an ordinary differential equation given in a compact precise form.  A further contribution of the paper is the proof that parameter convergence is guaranteed with the extremely weak interval excitation requirement.
\end{abstract}


\end{frontmatter}
%===============================================================================

\section{Introduction}
To comply with the stringent monitoring and control requirements in modern applications an accurate model of the system is vital. It is well-known that {\em nonlinear parameterizations} (NLP) are inevitable in any realistic dynamic model of practical problems with complex dynamics. Constitutive relations and
conservation equations used to characterize physical variables always involve NLP. Classical examples are friction dynamics \citep{ARMDUPDEW}, biochemical processes \citep{DOC} and in more recent technological developments we can mention fuel cells \citep{PUKSTEPENbook}, photovoltaic arrays \citep{MASbook}, windmill generators \citep{HEIbook} and biomechanics \citep{WINbook}. However, one of the assumptions that pervades almost all results in
adaptive estimation and control is {\em linearity} in the unknown parameters and there are very few results available for NLP systems. Quite often, in practical problems, there are only few physical parameters that are uncertain and occur nonlinearly in the underlying dynamic model. In some cases, it is possible to use suitable transformations so as to convert it into a problem where the unknown parameters occur linearly, usually involving overparameterizations. This procedure, however, suffers from serious drawbacks including the enlarging of dimension of the parameter space, with the subsequent increase in the excitation requirements needed to ensure parameter convergence. The reader is referred to \citep{ORTetalaut21} for a thorough discussion on the drawbacks of overparameterization.  

Some results for {gradient estimators} have been reported in the literature for {\em convexly} parameterized systems. It was first reported in \citep{FOMFRAYAKbook} (see also \citep{ORT}) that convexity is enough to ensure that the gradient search ``goes in the right direction" in a {\em certain region} of the estimated parameter space. The idea is then to apply a standard adaptive scheme in this region, while in the ``bad" region either the adaptation is frozen and a robust constant parameter controller is switched-on \citep{FRAetal} or, as proposed in \citep{ANNSKALOH}, the adaptation is running all the time and stability is ensured with a high-gain mechanism which is suitably adjusted incorporating prior knowledge on the parameters. In \citep{NETetal} {\em reparametrization} to convexify an otherwise non-convexly parameterized system is proposed. See also \citep{TYUetal,TYUetal1} for some interesting results along these lines, where the controller and the estimator switch between over/underbounding convex/concave functions. Some calculations invoking computationally demanding set membership principles---similar to fuzzy systems---have recently been reported in \citep{ADEGUALEH}. 

Using the Immersion and Invariance adaptation laws proposed in \citep{ASTKARORTbook}, stronger results were obtained in \citep{LIUetaltac,LIUetalscl} invoking the property of {\em monotonicity}, see also \citep{TYUetal,TYUetal1} for related results. The main advantage of using monotonicity, instead of convexity, is that in the former case the parameter search ``goes in the right direction" in {\em all regions} of the estimated parameter space---this is in contrast to the convexity-based designs where, as pointed out above, this only happens in some regions of this space. See the recent work \citep{ORTROMARA} where these results relying on monotonicity have been significantly extended. The reader is referred to \citep{ORTetalaut21,ORTROMARA} for recent reviews of the literature on parameter estimation and adaptive control of NLP systems.
Unfortunately, the monotonicity property can be exploited only for the case of {\em separable} NLP. That is for the case where we can factor the parameter dependent terms as $h_i(u,y,\theta)=\bar h_i(u,y)\psi_i(\theta_i)$, where $u$ and $y$ are measurable and $\theta_i$ is the unknown parameter. However, there are many practical application models where this factorization is not possible, we refer to this case as {\em non-separable} NLP. Two often encountered cases are $\cos(\theta_i \cdot h_i(u,y))$ or $e^{\theta_i \cdot h_i(u,y)}$. In particular, the last example appears in many physical processes including Arrenhius laws \citep{SILbook}, biochemical reactors \citep{DOC}, friction models \citep{ARMDUPDEW}, windmill systems \citep{BOBetalscl}, fuel cell systems \citep{XINetal}, photovoltaic arrays  \citep{BOBetalaut} and models of elastic moments \citep{SCHetal,SHAetal,YANDEQ}. This paper is devoted to the development of a systematic methodology for the parameter identification of systems containing this kind of exponential terms.
More precisely, we consider systems of the form 
\begalis{
	\dot x = F_x(u,y,\theta),
	y =H_x(u,y,\theta)
} 
with $u$ and $y$ measurable and $\theta$ a vector of unknown parameters, with some of its elements entering into the functions $F_x$ and/or $H_x$ via exponential terms of the form {\em $e^{\theta_i \cdot h_i(u,y)}$}. The objective is to design an on-line {\em estimator}
\begalis{
	\dot \chi = F_\chi(\chi,u,y),
	\hat \theta =H_\chi(\chi,u,y)
}
with $\chi(t) \in \rea^{n_\chi}$ such that we ensure {\em global exponential convergence} (GEC) of the estimated parameters. That is, for all $x(0) \in \rea^n,\chi(0) \in \rea^{n_\chi}$ and all continuous $u$ that generates a bounded state trajectory $x$ we ensure
\begequ
\lab{concon}
\liminf |\tilde \theta(t)| = 0,\quad \mbox{(exp)},
\endequ
where $\tilde \theta:=\hat \theta - \theta$ is the parameter estimation error, with all signals remaining bounded.

Notice that, in contrast with the existing approaches for non-separable NLP systems, we do not assume that the parameters live in known compact sets, that the nonlinearities satisfy some Lipschitzian properties, nor rely on injection of high-gain to dominate the nonlinearities or the use of complex, computationally demanding methodologies like $\mbox{min-max}$ optimizations, parameter projections or set membership techniques. Instead, we propose to design a classical on-line estimator whose dynamics is described by an ordinary differential equation given in a compact precise form. 

We identify in the paper two classes of systems for which the problem formulated above can be solved. The design procedure consists of the construction---from the non-separable NLP containing an exponential term---a new NLP regression equation (NLPRE) of the form $Y(u,y)=\phi^\top(u,y)\calg(\theta)$, where the functions $Y(u,y)$ and $\phi(u,y)$ are known and $\calg(\theta)$ is a nonlinear mapping. To estimate the parameters $\theta$ from the NLPRE we invoke the recent result of \citep{ORTROMARA}, where a least-squares plus dynamic regression equation and mixing \citep{ARAetaltac} (LS+DREM) estimator applicable for this kind of NLPRE is reported. A key feature of the LS+DREM estimator is that it ensures GEC imposing an extremely weak {\em interval excitation} (IE) assumption \citep{KRERIE,TAObook} of the regressor  $\phi$. On the other hand, this estimator requires that the mapping of the NLPRE satisfies a rather weak monotonizability property---that is captured by the verifiability of a linear matrix inequality (LMI) imposed on $\calg(\theta)$. We give two practical examples of the application of the proposed estimation method and illustrate their performance with some simulations.

\noindent {\bf Notation.} $I_n$ is the $n \times n$ identity matrix and $0_{s \times r}$ is an $s \times r$ matrix of zeros. $\rea_+$ and $\intnum_+$ denote the positive real and integer numbers, respectively. For $q \in \intnum_+$ we define the set $\bar q:=\{1,2,\dots,q\}$. For $a \in \rea^n$, we denote $|a|^2:=a^\top a$, and for any matrix $A$ its induced norm is $\|A\|$. All functions and mappings are assumed smooth and all dynamical systems are assumed to be forward complete. Given a function $h:  \rea^n \times \rea^m \to \rea$ we define its transposed gradient via the differential operator $\nabla_{(\cdot)} h(x,u):=\Big[\frac{\displaystyle \partial h }{\displaystyle \partial (\cdot)}(x,u)\Big]^\top$. For a mapping $\calg:\rea^{n_\eta} \to \rea^{p_\eta}$ we denote its Jacobian by $\nabla \calg(\eta):={\partial \calg \over \partial \eta}(\eta)$. To simplify the notation, the arguments of all functions and mappings are written only when they are first defined and are omitted in the sequel.

 %
%%%%%%%%%%%%%%%%%%5
\section{First Class of Systems}
\lab{sec2}
%%%%%%%%%%%%%%%%%%

In this section we consider NLP systems of the form
\begin{footnotesize}
\begsubequ
\lab{sys1}
\begali{
	\dot{x}&=f_1\left(x,u\right)+ f_2\left(x,u\right)\calg(\eta)
	\label{dotx1}\\
	y&=\begmat{y_1 \\ x}=\begmat{h_1\left(x,u\right)+ h_2\left(x,u\right)\theta_{2}+ h_3\left(x,u\right)e^{h_4\left(x\right)\theta_1}\\x} 
	\label{y1}
}
\endsubequ
\end{footnotesize}

with $x(t)\in {\mathbb{R}^n}$, $y(t)\in {\mathbb{R}^{n+1}}$ and $u(t)\in {\mathbb{R}^m}$ the systems state, output and control, respectively. The functions $f_i,\;i=1,2$, and $h_i,\;i=1,\dots,4$, are {\em known} nonlinear functions, $\calg:\rea^{n_\eta}\mapsto \rea^{p_\eta},\;p_\eta>n_\eta$, is a  {\em known} mapping of the unknown parameters $\eta \in \rea^{n_\eta}$,  and $\theta_{i} \in \rea,\;i=1,2$ are also unknown parameters. Hence, the overall vector of {\em unknown} parameters, which needs to be estimated on-line, consists of $\theta:=\col(\theta_1,\theta_2,\eta) \in \rea^{\ell_{I}}$, where $\ell_{I}:=2+n_\eta$. 

We make the important observation that, in view of the presence of the {\em exponential term} in the signal $y_1$, the parameterization of the system is  nonlinear and {\em non-separable}. As discussed in the Introduction none of the existing parameter estimators can deal with this difficult---but often encountered in practice---scenario. 
%
%%%%%%%%%%%%%%%5
\subsection{Assumptions}
\lab{subsec21}
We make the following assumptions on the system.

\begenu
\item[{\bf A1}]  [Sign definiteness] The scalar function $h_3$ is bounded away from {\em  zero}. That is $|h_3|>0$.

\item[{\bf A2}]  [Monotonicity] There exists a matrix $T_\calg \in \rea^{n_\eta \times p_\eta}$ such that the mapping $\calg(\eta)$ satisfies the {\em LMI}
\begequ
\lab{demcon1}
T_\calg \nabla \calg(\eta) + [\nabla \calg(\eta)]^\top T_\calg^\top \geq \rho_\calg I_{n_\eta},
\end{equation} 
for some $\rho_\calg {>0}$.
\endenu
%
%%%%%%%%%%%%%%%5
\subsubsection*{Discussion on the assumptions}
\noindent {\bf D1} 
In \citep[Proposition 1]{ORTetalaut21} it is shown that \eqref{demcon1} ensures  the mapping $T_\eta \calg(\eta) $ is {\em strictly monotonic} \citep{PAVetal}. That is, it satisfies
\begin{footnotesize}
	\begequ
\lab{monpro}
(a-b)^\top \left[T_\eta\calg(a) -T_\eta \calg(b)\right] \geq \rho_\eta |a-b|^2,\;\forall a,b \in \rea^{n_\eta},\; a \neq b.
\endequ 
\end{footnotesize}

This is the fundamental property that is required by the LS+DREM estimator used in the next section.

\noindent {\bf D2} 
The assumption that the state trajectories of \eqref{sys1} are bounded is standard in parameter estimation theory \citep{LJUbook,SASBODbook}. Similarly, the assumption that the dimension $n_\eta$ of the unknown parameters vector $\eta$ is smaller than $p_\eta$ is reasonable, otherwise we could redefine a new vector of unknown parameters $\bar \eta:=\calg(\eta) \in \rea^{n_\eta}$ without overparameterization and get a LRE. 
%
%%%%%%%%%%%%
\subsection{Regression Equation for Parameter Estimation}
\lab{subsec22}
In this section derive the regression equation that will be used to estimate the unknown parameters $\theta$. As expected, this regressor equation is nonlinearly parameterized, which hampers the application of standard estimation techniques. Therefore, we are compelled to appeal---in Section \ref{sec4}---to the  LS+DREM parameter estimator recently reported in \citep{ORTROMARA,PYRetal}. 

\begin{lemma}\em
	\lab{lem1}
	Consider the system \eqref{sys1} verifying  Assumptions {\bf A1, A2}. There exists measurable, scalar signals $Y_I(x,u,y),\;\phi_{I,i}(x,u,y),\;i=1,\dots,s_I,\;s_I:=3+2p_\eta,$ such that the following NLPRE holds:
	\begin{equation}
		\label{nlpre1}
		Y_I(x,u,y)= \phi_I^\top(x,u,y) \calw_I(\theta),
	\end{equation}
	where we defined the mapping $\calw_I:\rea^{\ell_{I}} \to \rea^{s_I}$
	\begequ
	\lab{calw1}
	\calw_I(\theta):=\begmat{\theta_1 & \theta_2 & \theta_1\theta_2 & \theta_1 \calg^\top(\eta) & \theta_1\theta_2\calg^\top(\eta)}^\top.
	\endequ
\end{lemma}

%%%%%%%%%%%%%%%5
\subsubsection*{Discussion on the regressor equation}
\noindent {\bf D3} 
It is possible to construct another NLPRE proceeding as follows. First, exploiting the monotonicity property of Assumption {\bf A2} and using the LS+DREM algorithm estimate the parameters $\eta$ filtering \eqref{dotx1}. Then, use this estimate in the (approximate) calculation of $\dot h_4$, yielding
$$
\dot {\hat h}_4 = \nabla^\top h_4 [f_1 + f_2 \calg(\hat \eta)].
$$
Applying the certainty equivalent principle, and replacing this expression in the chain of implications of the proof of Lemma \ref{lem1} in Appendix A would then yield a simpler NLPRE where only the terms $(\theta_1,\theta_2, \theta_1 \theta_2)$ will appear. Of course, the drawback of this approach is that we rely on the fast convergence of $\tilde \eta:=\hat \eta - \eta$ to zero. 

\noindent {\bf D4}   
In the system \eqref{sys1} the function $h_4$ appearing in the exponential does not depend on $u$. It is possible to adapt the result of Lemma \ref{lem1} to consider that case in the following way. The expression for $\dot h_4$ given in \eqref{doth4} would need to be replaced by
$$
\dot h_4 = \nabla_x^\top h_4 [f_1 + f_2 \calg(\eta)]+\nabla_u^\top h_4 \dot u.
$$
To construct the NPLRE as in Lemma \ref{lem1} for this case it is clearly necessary to know $\dot u$. However, in many practical applications the control law contains an {\em integral action}---{\em e.g.}, in PID control---therefore this signal is available for measurement. 
%
%%%%%%%%%%%%
\subsection{Construction of a Strictly Monotonic Mapping}
\lab{subsec23}
%%%%%%%%%%%
%
To estimate the parameters $\theta$ from the NLPRE    \eqref{nlpre1} we invoke the recent result of \citep{ORTROMARA}, where the LS+DREM estimator proposed in \citep{PYRetal}, which is applicable for linear regression equations,  was {\em extended} to deal with NLPRE. However, this estimator requires that the mapping of the NLPRE satisfies a {\em monotonicity} property, which is not verified by $\calw_I(\theta)$ given in \eqref{calw1}. Therefore, in this section we construct a new mapping verifying the required monotonicity condition.

\begin{lemma}
	\lab{lem2}\em
	Consider the mapping $\calw(\theta)$ given in  \eqref{calw1} with $\calg(\eta)$ verifying Assumption {\bf A2}. There exists a constant $\alpha_m>0$ such that for all $\alpha \geq \alpha_m$ the mapping $\calw_I(\theta)$ satisfies the {\em LMI}
	\begequ
	\lab{demcon1}
	T_{\calw_I}\nabla \calw_I(\theta) + [\nabla \calw_I(\theta)]^\top T_{\calw_I}^\top \geq \rho_{\calw_I} I_{\ell_{I}},
\end{equation}
for some $\rho_{\calw_I}>0$, with the matrix 
$$
T_{\calw_I}:=\begmat{\alpha & 0 & 0 &  0_{1 \times p_\eta}  &  0_{1 \times p_\eta}\\ 0 & \alpha & 0 &  0_{1 \times p_\eta}  &  0_{1 \times p_\eta}\\ & 0_{n_\eta \times 3} & & \emph{sign}(\theta_1) T_\calg &  0_{n_\eta \times p_\eta}} \in \rea^{\ell_{I}\times s_I}.
$$
\end{lemma}
%
%%%%%%%%%%%%%%%5
\subsubsection*{Discussion on the mapping}
\noindent {\bf D5} 
Notice that the only prior knowledge needed to construct the matrix $T_{\calw_I}$ is $\sign(\theta_1)$. On the other hand, to select the value of $\alpha$  some prior knowledge on the parameters $\theta$  is required. Specifically, as shown in the proof of Lemma \ref{lem2} in Appendix A,  it is necessary to know an upper bound on $\|  T_\calg \calg(\eta)\|$.
%%
%%%%%%%%%%%%%%%%%%5
\section{Second Class of Systems}
\lab{sec3}
%%%%%%%%%%%%%%%%%%
%
In this section we consider second order systems of the form
\begsubequ
\lab{sys2}
\begali{
	\lab{dotx2}
	\ddot x & = f_1(x)+f^\top_2(x,\dot x)\calg(\eta) + h_3(x) e^{\theta_1 h_4(x)}+ u\\
	\lab{y2}
	y&=\begmat{x \\ \dot x}
}
\endsubequ
with $x(t)\in {\mathbb{R}}$ and $u(t)\in {\mathbb{R}}$. The functions $f_i,\;i=1,2$, and $h_i,\;i=1,3$, are {\em known} nonlinear functions, $\calg:\rea^{n_\eta}\mapsto \rea^{p_\eta},\;p_\eta>n_\eta$, is a  {\em known} mapping of the unknown parameters $\eta \in \rea^{n_\eta}$,  and $\theta_{1} \in \rea$ is also an unknown parameter. Hence, the overall vector of unknown parameters, which needs to be estimated on-line, consists of $\theta:=\col(\theta_1,\eta) \in \rea^{\ell_{II}}$, where $\ell_{II}:=1+n_\eta$. 

Notice that, in contrast to system \eqref{sys1}, in this case the dynamics is second order and the nasty exponential term enters into the state question instead of the readout map. Moreover, note that the control signal is scalar and enters linearly in the state equation. In particular, observe that the function $h_3$ appearing in the exponential {\em does not} depend on $u$ now.\footnote{To simplify the presentation, but with an obvious abuse of notation, we keep the same symbol for both functions.} 

To simplify the calculations, in the model \eqref{sys2} we do not include unknown parameters multiplying the function $h_3$ or the control $u$. As explained in Discussion {\bf D7} below, this can be easily added redefining $h_3(x):=\theta_2 \bar h_3(x)$  and  $u:=\theta_3 \bar u$, where the functions $\bar h_3$ and $\bar u$ are known but $\theta_2$ and $\theta_3$ are unknown parameters.  
%
%%%%%%%%%%%%%%%5
\subsection{Assumptions}
\lab{subsec31}
We make  on the system \eqref{sys2} Assumptions {\bf A1, A2} together with the following.

\begenu
\item[{\bf A3}]  [Separability] The function $f_2(x,\dot x)$ verifies
$$
\nabla_{\dot x}f_2(x,\dot x)=\psi_a(x)\psi_b(\dot x),
$$
for some functions $\psi_a(x)$ and $\psi_b(\dot x)$.
\endenu
%
%%%%%%%%%%%%%%%5
\subsubsection*{Discussion on Assumption {\bf A3}}
\noindent {\bf D6} 
As shown in the proof of Lemma \ref{lem3} given in Appendix A, Assumption {\bf A3} is needed to be able to generate---via LTI filtering---a measurable regressor in the NLPRE. We observe that the function $\nabla_{\dot x}f_2 \in \rea^{p_\eta}$ hence, for $p_\eta >1$, this is a vector function. However, there is no restriction on the dimensions of the functions $\psi_a$ and  $\psi_b$, as long as they comply with the dimensionality requirement $\psi_a\psi_b \in \rea^{p_\eta}$. This degree of freedom relaxes the condition of the assumption. 
%%%%%%%%%%%%
\subsection{Regression Equation for Parameter Estimation}
\lab{subsec32}
As in Subsection \ref{subsec22} we derive here the NLPRE that will be used to estimate the unknown parameters $\theta$. 

\begin{lemma}\em
	\lab{lem3}
	Consider the system \eqref{sys2} verifying  Assumptions {\bf A1-A3}. There exists measurable, scalar signals $Y_{II}(x,u,y),\;\phi_{II,i}(x,u,y),\;i=1,\dots,s_{II},\;s_{II}:=1+2p_\eta,$ such that the following NLPRE holds:
	\begin{equation}
		\label{nlpre2}
		Y_{II}(x,u,y)= \phi_{II}^\top(x,u,y) \calw_{II}(\theta),
	\end{equation}
	where we defined the mapping $\calw_{II}:\rea^{\ell_{II}} \to \rea^{s_{II}}$
	\begequ
	\lab{calw2}
	\calw_{II}(\theta):=\begmat{\theta_1 &  \calg^\top(\eta) & \theta_1\calg^\top(\eta)}^\top.
	\endequ
\end{lemma}
%
%%%%%%%%%%%%%%%5
\subsubsection*{Discussion on regression equation}
\noindent {\bf D7} 
To include an unknown multiplicative parameter in the function $h_3$ or the control $u$ we proceed as follows. Define $h_3(x)=\theta_2 \bar h_3(x)$ and $u=\theta_3 \bar u$, where the functions $\bar h_3$ and $\bar u$ are known but $\theta_2$ and $\theta_3$ are unknown parameters.  Tracing back  the proof of Lemma \ref{lem3} given in Appendix A, in the first step where we divide the model equation by $h_3$ we divide instead by $\bar h_3$. Then, the parameter $\theta_2$ appears multiplying the exponential in the term in parenthesis and it is removed in the next line. That is, the first three lines of the proof become
\begali{
	\nonumber
	{1 \over \bar h_3} \ddot x  = \theta_2 e^{h_4\theta_1}  + \bar f_3+\bar f^\top_4 \calg(\eta) +  \theta_3 {\bar u \over \bar h_3}}
\begali{
	\nonumber
	\dif - {\dot {\bar h}_3 \over \bar h^2_3}  \ddot x+  {1 \over \bar h_3}{d^3 x \over dt^3} =  \theta_1 \dot h_4 \Big(\theta_2 e^{h_4\theta_1}\Big) + \dot {\bar f}_3 + \dot {\bar f}^\top _4 \calg(\eta) -}
\begali{
	\nonumber
	  \theta_3 \Big( {\dot {\bar h}_3 \over \bar h^2_3} \bar u - {\dot {\bar u} \over \bar h_3} \Big)}
\begali{ 
	\nonumber
	\Longleftrightarrow - {\dot {\bar h}_3 \over \bar h^2_3}  \ddot x+  {1 \over \bar h_3}{d^3 x \over dt^3} =  \theta_1 \dot h_4 \Big( {1 \over \bar h_3} \ddot x - \bar f_3 - \bar f^\top_4 \calg(\eta) -\theta_3 {\bar u \over \bar h_3} \Big) }
\begali{
	\nonumber
	+\dot {\bar f}_3 + 	 \dot {\bar f}^\top _4 \calg(\eta) -  \theta_3 \Big( {\dot {\bar h}_3 \over \bar h^2_3} \bar u - {\dot {\bar u} \over \bar h_3} \Big),  
}
with the new definitions  
$$
\bar f_3:= {f_1 \over \bar h_3} ,\; \bar f_4:= {1 \over \bar h_3}f_2.
$$ 
The remaining part of the proof remains unchanged leading to a NLPRE similar to \eqref{nlpre2}, with the new $\bar {(\cdot)}$ terms and adding to the parameter vector $\theta_3$ and $\theta_1\theta_3$. As proven in Proposition \ref{pro1}, from this NLPRE we can estimate exponentially fast $(\theta_1,\theta_3,\eta)$. Therefore, we can replace their estimates in the model \eqref{sys2} leading to  the system
$$
\ddot z  = f_1(z)+f^\top_2(z,\dot z)\calg(\hat \eta) + \theta_2 \bar h_3(z) e^{\hat \theta_1 h_4(z)}+ \hat \theta_3 \bar u,
$$
which is a classical linearly parameterized system from which we can estimate $\theta_2$ with standard filtering plus gradient descent techniques.
%
%%%%%%%%%%%%

\subsection{Construction of a Strictly Monotonic Mapping}
\lab{subsec33}
%%%%%%%%%%%
%
Similarly to the calculations presented in Subsection \ref{subsec23} we present here the matrix $T_{\calw_{II}} \in \rea^{\call_{II} \times s_{II}}$ that defines the new monotonic mapping. The proof of this lemma is trivial, therefore it is omitted for brevity.

\begin{lemma}
	\lab{lem4}\em
	Consider the mapping $\calw_{II}(\theta)$ given in  \eqref{calw2} with $\calg(\eta)$ verifying Assumption {\bf A2}. The mapping $\calw_{II}(\theta)$ satisfies the {\em LMI}
	\begequ
	\lab{demcon2}
	T_{\calw_{II}}\nabla \calw_{II}(\theta) + [\nabla \calw_{II}(\theta)]^\top T_{\calw_{II}}^\top \geq \rho_{\calw_{II}} I_{\ell_{II}},
\end{equation}
with the matrix
$$
T_{\calw_{II}}:=\begmat{1 &  0_{1 \times p_\eta}  &  0_{1 \times p_\eta}\\ 0_{n_\eta \times 1} &   T_\calg &  0_{n_\eta \times p_\eta}} \in \rea^{\ell_{II}\times s_{II}}.
$$
\end{lemma}
\subsubsection*{Discussion on the mapping}
\noindent {\bf D8} 
Notice that, in contrast with the construction of Subsection \ref{subsec23}, here there is no requirement of prior knowledge on the parameter $\theta_1$. These stems from the fact that, as seen in  \eqref{calw2}, the mapping $\calg(\eta)$ appears once without multiplying this parameter---compare with \eqref{calw1}. Therefore, Assumption {\bf A2} is sufficient to construct the new monotonic mapping. 
%%

%%%%%%%%%%%%

\section{A Globally Exponentially Convergent  Estimator of $\theta$}
\lab{sec4}
%%%%%%%%%%%

In this section we present the main result of the paper, that is, an estimator of the parameters $\theta$ that achieves GEC of the parameter error.  We proceed from the NLPREs constructed in Lemmata \ref{lem1} and \ref{lem3} and, as explained in Subsection \ref{subsec23}, we propose to use the  LS+DREM estimator recently reported in \citep{ORTROMARA}. Towards this end, we use the new mappings identified in  Lemmata \ref{lem2} and \ref{lem4} that verify the monotonicty conditions required by the LS+DREM estimator. To simplify the notation we avoid the subindices $(\cdot)_{I}$ and  $(\cdot)_{II}$ of the various terms appearing in previous sections and present a single proposition applicable to both classes of systems.

Therefore, we consider a general scalar NLPRE of the form
\begequ
\lab{nlpre}
Y(t)= \phi^\top(t) \calw(\theta)
\endequ
with $\calw: \rea^\ell \to \rea^s$. The main feature of the LS+DREM estimator is that it ensures GEC imposing the following extremely weak {\em IE} assumption \citep{KRERIE,TAObook} of the regressor  $\phi$.\\

\noindent {\bf A4} [Excitation] The regressor vector $\phi$ is IE. That is, {there exist} constants $C_c>0$ and $t_c>0$ such that
$$
\int_0^{t_c} \phi(s) \phi^\top(s)  ds \ge C_c I_s.
$$
The proof of the proposition below is given in \citep[Proposition 1]{ORTROMARA}, therefore it is omitted here. 

\begin{proposition}\em
	\lab{pro1}
	Consider the NLPRE \eqref{nlpre} with $\phi$ verifying  Assumption {\bf A4} and $\calw$ satisfying the LMI
	$$
	T_{\calw}\nabla \calw(\theta) + [\nabla \calw(\theta)]^\top T_{\calw}^\top \geq \rho_{\calw} I_{\ell}
	$$
	for some matrix $T_{\calw} \in \rea^{\ell \times s}$ and $\rho_{\calw}>0$. Define the LS+DREM interlaced estimator
	\begalis{
		\dot{\hat \calw} & =\gamma_\calw F \phi (Y-\phi^\top \hat\calw),\; \hat\calw(0)=\calw_{0} \in \rea^{s}\\
		\dot {F}& =  -\gamma_\calw F \phi  \phi^\top  F,\; F(0)={1 \over f_0} I_s \\
		\dot{\hat \theta} & = \Delta \Gamma T_\calw  [\caly -\Delta \calw(\hat\theta) ],\; \hat\theta(0)=\theta_{0} \in \rea^\ell,
	}	
	with tuning gains the scalars $\gamma_\calw>0$, $f_0>0$ and the positive definite matrix $\Gamma \in \rea^{\ell\times \ell}$,  and we defined the signals
	\begalis{
		\Delta :=\det\{I_s-f_0F\},
		\caly := \emph{\adj} \{I_s- f_0F\} (\hat\calw -  f_0F \calw_{0}),
	}
	where $ \emph{adj}\{\cdot\}$ denotes the adjugate matrix. For all initial conditions $\calw_{0} \in \rea^s$ and $\theta_{0} \in \rea^\ell$. The estimation errors of the parameters $\tilde \theta$ verify \eqref{concon} with all signals {\em bounded}. 
\end{proposition}
%
%
%%%%%%%%%%%%%%%5
\section{Two Practical Examples}
\lab{sec5}
%%%%%%%%%%%%
%
\subsection{Proton Exchange Membrane Fuel Cell}
\lab{subsec51}
%%%%%%%%%%%%
%
Parameter estimation is vital for modeling and control of fuel cell systems. However, an accurate description of the fuel cell dynamics implies the use of models with nonlinear parameterizations \citep{PUKSTEPENbook}.  The interested reader is refered to  \citep{XINetal} where a detailed review of the literature on parameter estimation of fuel cell systems is reported.

\subsubsection*{Verification of the conditions from the general result}
A widely accepted mathematical model of a {\em Proton Exchange Membrane Fuel Cell} (PEMFC) is given in \citep[Section II.B]{XINetal}. It can be shown that this model can be written in the form \eqref{sys1} with $n=m=n_\eta=p_\eta=1$ and the scalar linear map $\calg(\eta)=\eta$.  

We make the observation that  function $h_3$, which is proportional to the membrane conductivity, is bounded away from zero, hence verifying Assumption {\bf A1}.

Since $\calg=\theta_3$ the mapping $\calw_I:\rea^3 \to \rea^5$ defined in \eqref{calw1} is simpler and given as
$$
\calw_I(\theta):=\begmat{\theta_1 & \theta_2& \theta_1\theta_2 & \theta_1\theta_3 & \theta_1\theta_2\theta_3}^\top.
$$

Some simple calculations give us terms $Y_I$ and $\phi_I^\top$ for  the NLPRE \eqref{nlpre1}.
And the  matrix $T_{\calw_I}$ of Lemma \ref{lem2} is given as
$$
T_{\calw_I}:=\begmat{\alpha & 0 & 0 & 0 & 0 \\ 0 & \alpha & 0 & 0 & 0\\ 0 & 0 & 0 & \sign(\theta_1) & 0 },
$$
and the minimum value for $\alpha$ is $
\alpha_m =  {\theta^2_3 \over |\theta_1|}.$

	\subsection{Human Shank Dynamics}
	\lab{subsec52}
	%%%%%%%%%%%%
	Neuromuscular electrical stimulation is an active research area that aims at restoring functionality to human limbs with motor neuron disorders. Control of these systems is a challenging problem because the musculoskeletal dynamics are nonlinear and highly uncertain \citep{WINbook}. In this subsection we are  interested in the mechanical dynamics of the human shank motion where the input is the joint torque produced by electrode stimulation of the shank muscles. We consider the scenario described in detail in \citep{YANDEQ}, see also \citep{SCHetal, SHAetal} and concentrate our attention on the problem of estimating the parameters of a widely accepted mathematical model of this system. Namely, the system described by equations (11) to (14) of \citep{YANDEQ}, that we repeat here for ease of reference 
	\begin{footnotesize}
	\begequ
	\lab{jddotx}
	J \ddot x +b_1 \dot q + b_2 \sign(\dot x)+ k_1 e^{-k_2 x}(x-q_0)+mg \ell \sin(x)=u,
	\endequ
\end{footnotesize}

	where $(x,\dot x)$ are assumed measurable and all the parameters are assumed unknown. The reader is referred to this reference for further details on the model, in particular, the physical interpretation of the different terms in the model, and the overall formulation of the neuromuscular electrical stimulation problem.  
	
	\subsubsection*{Verification of the conditions from the general result}
	The following clarifications regarding our formulation of the parameter estimation problem are in order.
	\begite
	\item[{\bf C1}]
	As indicated in  \citep{YANDEQ},  the term $\sign(\dot x)$ of our model \eqref{jddotx} is replaced in equation (12) of  \citep{YANDEQ} by the function $\tanh(b_3\dot x)$, with a large value for $b_3>0$, which is a smooth {\em approximation} of the sign function. This approximation is made for mathematical convenience of their calculations that rely on a smoothness assumption, but is not required in our approach that can deal with discontinuous nonlinearities.
	
	\item[{\bf C2}] In this paper we assume that the term $q_0$, which is the constant resting knee angle, and the constant inertia $J$ are {\em known}.  Therefore the uncertain parameters in our case are $\col(b_1,b_2,k_1,k_2,m\ell)$. The assumption of known $J$ is not to restrictive because the inertia can be predicted
	from the subject's anthropometric data \citep{WINbook}.
	
	\item[{\bf C3}] In  \citep{YANDEQ} there is an additional, bounded, unstructured, additive term in \eqref{jddotx} that is omitted here for brevity. As shown in Proposition \ref{pro1} we achieve GEC of the parameter estimates, therefore this term could be easily accommodated in our analysis to ensure practical stability. 
	\endite
	
	The dynamics \eqref{jddotx} belongs to the second class of systems given by \eqref{sys2} with $n_\eta=p_\eta=3$, and the following definitions for the functions
	\begalis{
		& f_1(x)=0,\;f_2(x,\dot x)={1 \over J}\col(-\dot x,-\sign(\dot x),g\sin(x)),\; \\
		& h_3(x)=k_1 \bar h_3(x):=k_1 {1 \over J}(x-q_0),\;h_4(x)=-x,
	} 
	and the parameters
	\begalis{
		& \theta_1=k_2,\; \calg(\eta)=\eta=\col(b_1,b_2,m\ell),\;\theta:=\col(\theta_1,\eta^\top).
	}
	We bring to the readers attention the fact that the model \eqref{jddotx} has a parameter $k_1$ multiplying the exponential term. Therefore, it is necessary to invoke the two-stage certainty-equivalent based procedure described in Discussion {\bf D7}. That is, we estimate with the NLPRE \eqref{nlpre2} the parameters $(k_2,b_1,b_2,m\ell)$ and then estimate, {\em e.g.}, with some filtering and a standard gradient, the remaining parameter $k_1$.
	
	To comply with Assumption {\bf A1}, we assume that $|x-q_0|>0$.\footnote{Adding a simple logic and a discontinuous function we can easily avoid the singularity points and replace this assumption by the knowledge of a set such that $q_0 \in [q_0^m,q_0^M]$.} Clearly, since $\calg(\eta)=\eta$, Assumption {\bf A2} is satisfied with $T_{\calg}={\rho_\calg \over 2} I_3$, with any $\rho_\calg>0$. Finally  Assumption {\bf A3} is satisfied with the functions
	$$
	\psi_a(x):= \begmat{-1 & 0 & 0 \\ 0 & -1 & 0 \\ 0 & 0 & g \sin(x)},\;\psi_b(\dot x):= \begmat{\dot x \\ \sign(\dot x) \\ 1}.
	$$
	
	The mapping $\calw_{II}:\rea^4 \to \rea^7$ is  given as
	$$
	\calw_{II}(\theta):=\begmat{k_2 & b_1  & b_2 & m\ell & k_2 b_1  & k_2 b_2 & k_2 m\ell}^\top.
	$$
	
	Some simple calculations give us terms $Y_{II}$ and $\phi_{II}^\top$ for  the  NLPRE \eqref{nlpre2}.
	
	Finally, the  matrix $T_{\calw_{II}} \in \rea^{4 \times 7}$ of Lemma \ref{lem4} is given as
	$$
	T_{\calw_{II}}:=\begmat{1 & 0 & 0 & 0 & 0 & 0 & 0 \\ 0 & {\rho_\calg \over 2} & 0 & 0 & 0 & 0 & 0 \\ 0 & 0 & {\rho_\calg \over 2} & 0 & 0 & 0 & 0 \\ 0 & 0 & 0 & {\rho_\calg \over 2} & 0 & 0 & 0}.
	$$
	
		\section{Simulation Results}

\subsection{First class of systems}
\lab{subsec51}
%%%%%%%%%%%%

Consider the ''synthetic'' model of the first slass of the systems in the form \eqref{sys1}, where $f_1=-x+u$,  $f_2=u$, $h_1=x^2$, $h_2=x-353$, $h_3= 0.1+x^2$, $h4=-1/x$.

Since $\calg=\theta_3$ the mapping $\calw_I:\rea^3 \to \rea^5$ defined in \eqref{calw1} is simpler and given as
$$
\calw_I(\theta):=\begmat{\theta_1 & \theta_2& \theta_1\theta_2 & \theta_1\theta_3 & \theta_1\theta_2\theta_3}^\top.
$$

Some simple calculations show that the terms of the NLPRE \eqref{nlpre1} are given as
\begalis{
	Y_I&= F(p)p\Big({y-h_1 \over h_3}\Big)\\
	\phi_I^\top&=F(p)\begmat{p\Big({h_2 \over h_3}\Big) &{f_1 (y-h_1) \over x^2h_3}&{f_2 (y-h_1) \over x^2h_3}&{f_1h_2 \over x^2h_3}&-{f_2h_2 \over x^2h_3}}.
}

For simulations we used next parameter values: filters parameter $\lambda=600$, $\hat{\calw}(0)= [0 \; 0 \; 0 \; 0 \; 0]$, $f_0=1$,
$\hat{\theta}=col[0 \; 0 \; 0]$, $\Gamma =10^{3} \begin{bmatrix}
	&1 &0 &0 \\
	&0 &1 &0 \\
	&0 &0 &1	
\end{bmatrix}$, 
$T_{\calw}=
\begin{bmatrix}
	&10 &0 &0 &0 &0\\
	&0 &10 &0 &0 &0\\
	&0 &0 &0 &10 &0
\end{bmatrix}$.

\begin{figure}[ht]
	\centering
	\includegraphics[width=1\linewidth]{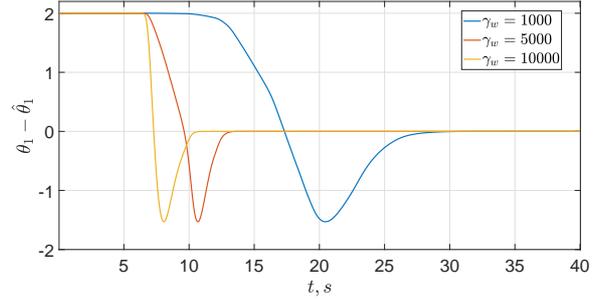}
	\caption{Transients of the error $\theta_1 -\hat{\theta}_1$ for different values of the parameter $\gamma_w$}
	\label{err_tet_1}
\end{figure}

\begin{figure}[ht]
	\centering
	\includegraphics[width=1\linewidth]{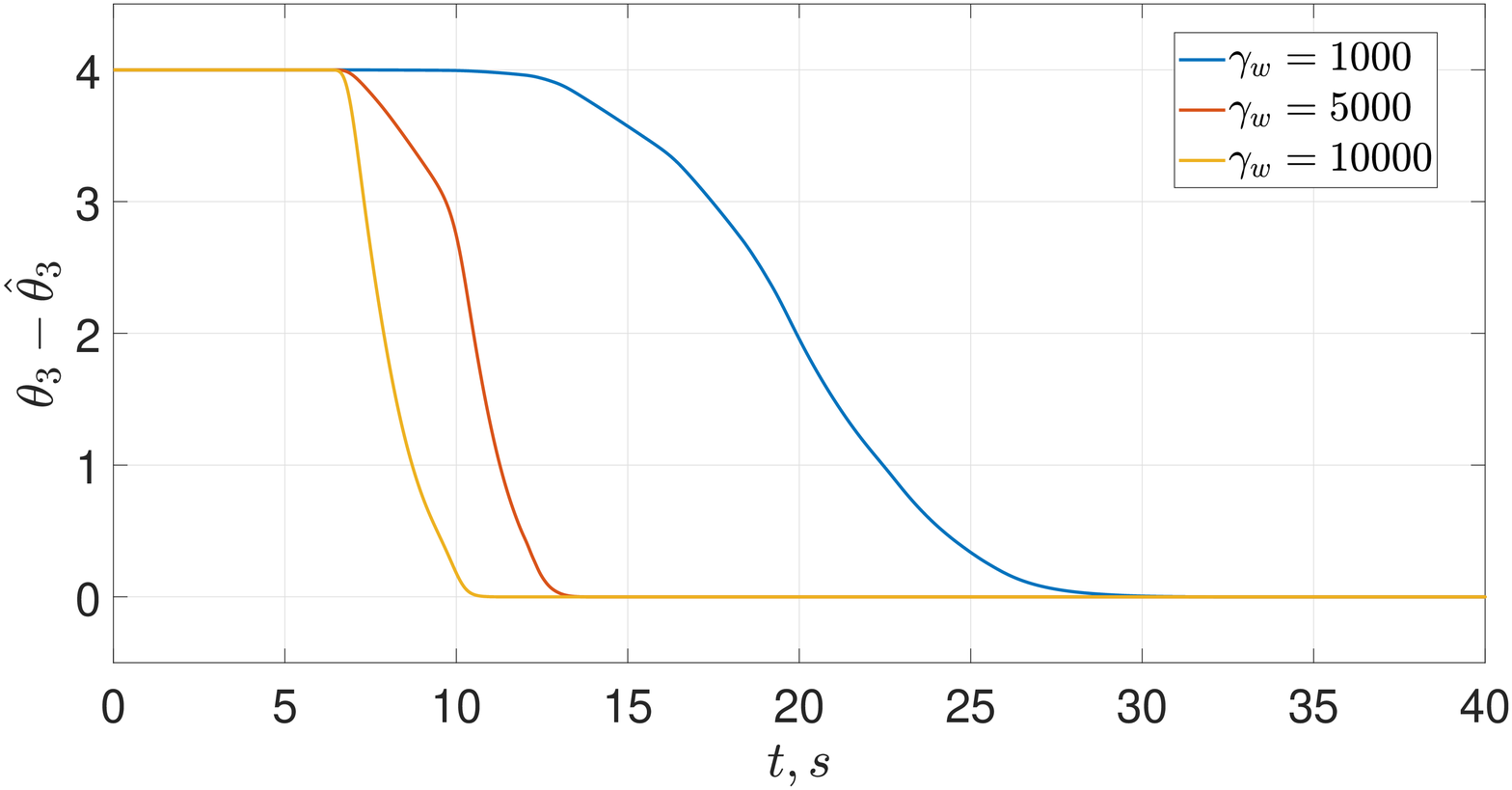}
	\caption{Transients of the error $\theta_3 -\hat{\theta}_3$ for different values of the parameter $\gamma_w$}
	\label{err_tet_3}
\end{figure}

Fig. \ref{err_tet_1} ... Fig. \ref{err_tet_3} demonstrate error transients of parameter estimations. In simulations we swithed on our observer on the fifth second. Figures demonstrate that  error transients  tends to the zero.

\begin{figure}[ht]
	\centering
	\includegraphics[width=1\linewidth]{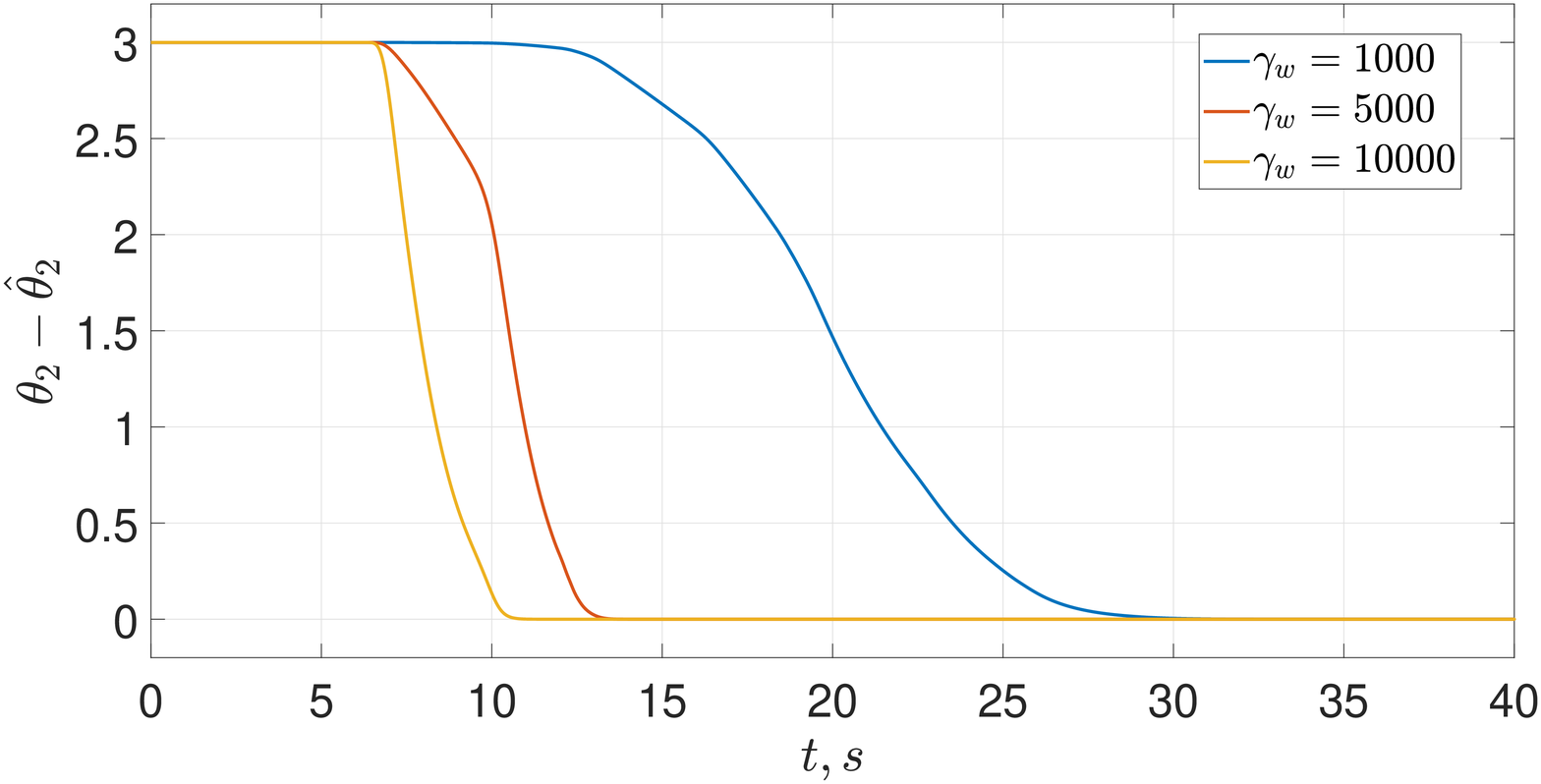}
	\caption{Transients of the error $\theta_2 -\hat{\theta}_2$ for different values of the parameter $\gamma_w$}
	\label{err_tet_2}
\end{figure}

	\subsection{Second class of systems}
	\lab{subsec52}
	%%%%%%%%%%%%
	
	Parameters of the human shank model were chosen as in \citep{YANDEQ}. For estimation of shank model parameters we used algorithm from proposition 1 with
	$\hat{\calw}(0)= [0 \; 0 \; 0 \; 0 \; 0 \; 0 \; 0]$, $f_0=10^{-3}$,
	$\Gamma = \begin{bmatrix}
		&1 &0 &0 &0\\
		&0 &1 &0 &0 \\
		&0 &0 &1 &0\\
		&0 &0 &0 &1
	\end{bmatrix}$, 
	$T_{\calw}=
	\begin{bmatrix}
		&1 &0 &0 &0 &0 &0 &0\\
		&0 &1 &0 &0 &0 &0 &0\\
		&0 &0 &1 &0 &0 &0 &0\\
		&0 &0 &0 &1 &0 &0 &0
	\end{bmatrix}$, 
	$\hat{\theta}=col[0 \; 0 \; 0 \; 0 \; 0 \; 0 \; 0]$ and $\lambda=1$ in the filters.
	
Consider computer simulations of human shank system with dynamic robust control from \citep{YANDEQ}
\begin{align}
	\nonumber
	u(t)=(k_1+1)r(t)-(k_1+1)r(0)+\\
	\nonumber
	+\int\limits_0^t \left[(k_1+1)k_2 r(\tau)+k_3 sign(r(\tau))\right]d\tau
\end{align}
where
$e=x-x$ and $r=\dot{e}-\mu e$.
For simulation we used
$\dot{x}_1=
\begin{bmatrix}
	&0 &1 &0\\
	&0 &0 &1\\
	&-6000 &-1300 &-80
\end{bmatrix}
x_1+
\begin{bmatrix}
	0\\0;\\1
\end{bmatrix}
w$, 
\begin{equation}
	\left\{
	\begin{matrix}
		&1000 \pi/3, &0\leqslant t < 20 s,\\
		&1000 \pi, &10 \leqslant t < 40 s,\\
		&2500 \pi, &t\geqslant40 s,
	\end{matrix}
	\right.
\end{equation}

$x_1=[x_d\;\dot{x}_d \; \ddot{x}_d]^\top$.

\begin{figure}[ht]
	\centering
	\includegraphics[width=1\linewidth]{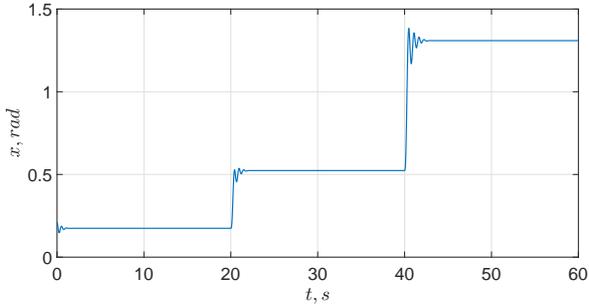}
	\caption{Transient of the output $x$}
	\label{x}
\end{figure}

\begin{figure}[ht]
	\centering
	\includegraphics[width=1\linewidth]{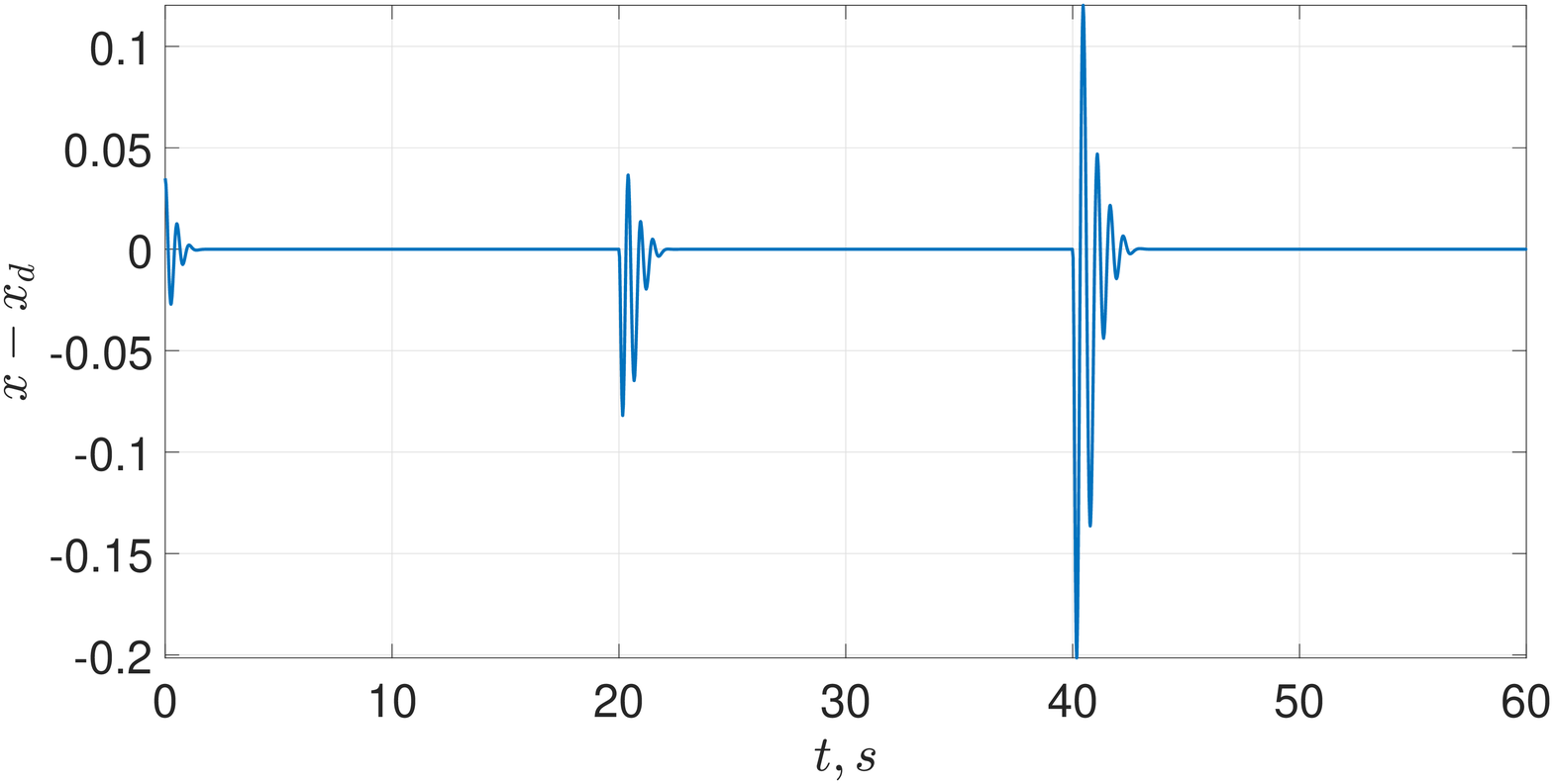}
	\caption{Transients of the error $x-x_d$}
	\label{err}
\end{figure}

For simulation we used $\mu=4$, $k_1=1$, $k_2=2$ and $k_3=40$ for control algorithm (simulation results for Human Shank with control are shown on Fig. \ref{x} and Fig. \ref{err}) and $\lambda=10$, $\gamma=10^6$, $\Gamma=10^5I$, $f_0=0.001$. 

\begin{figure}[ht]
	\centering
	\includegraphics[width=1\linewidth]{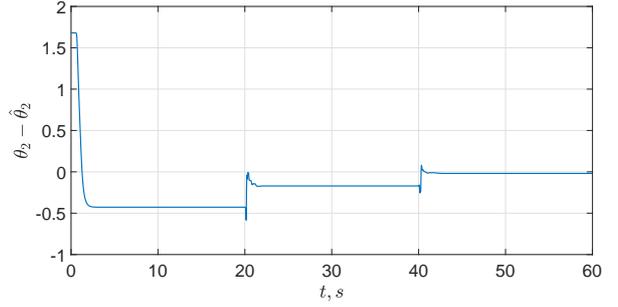}
	\caption{Transients of the error $\theta_2-\hat{\theta}_2$ ($k_2-\hat{k}_2$)}
	\label{2_hat_eta_1}
\end{figure}

\begin{figure}[ht]
	\centering
	\includegraphics[width=1\linewidth]{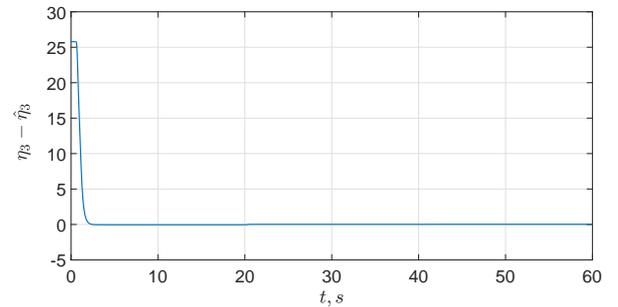}
	\caption{Transients of the error $\eta_3-\hat{\eta}_3$ ($\frac{m g l}{J}-\widehat{\frac{m g l}{J}}$)}
	\label{2_hat_eta_4}
\end{figure}

Fig. \ref{x} demonstrates transient of the system output $x$. Fig. \ref{err} demonstrates transient of the error $x-x_d$. Fig. \ref{2_hat_eta_1}...\ref{2_hat_eta_4} demonstrate transients of estimation errors.

\begin{figure}[ht]
	\centering
	\includegraphics[width=1\linewidth]{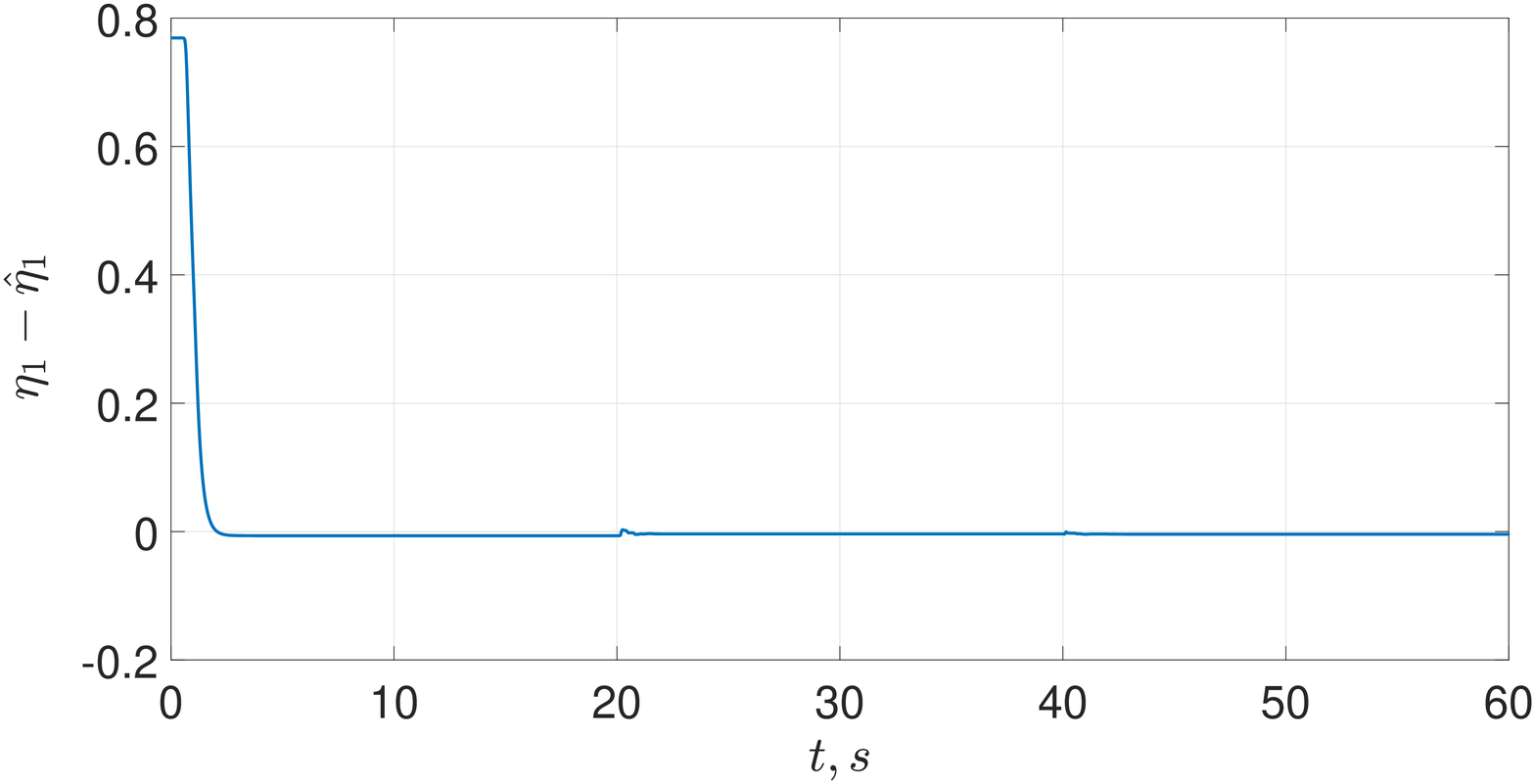}
	\caption{Transients of the error $\eta_1-\hat{\eta}_1$ ($\frac{b_1}{J}-\widehat{\frac{b_1}{J}}$)}
	\label{2_hat_eta_2}
\end{figure}

\begin{figure}[ht]
	\centering
	\includegraphics[width=1\linewidth]{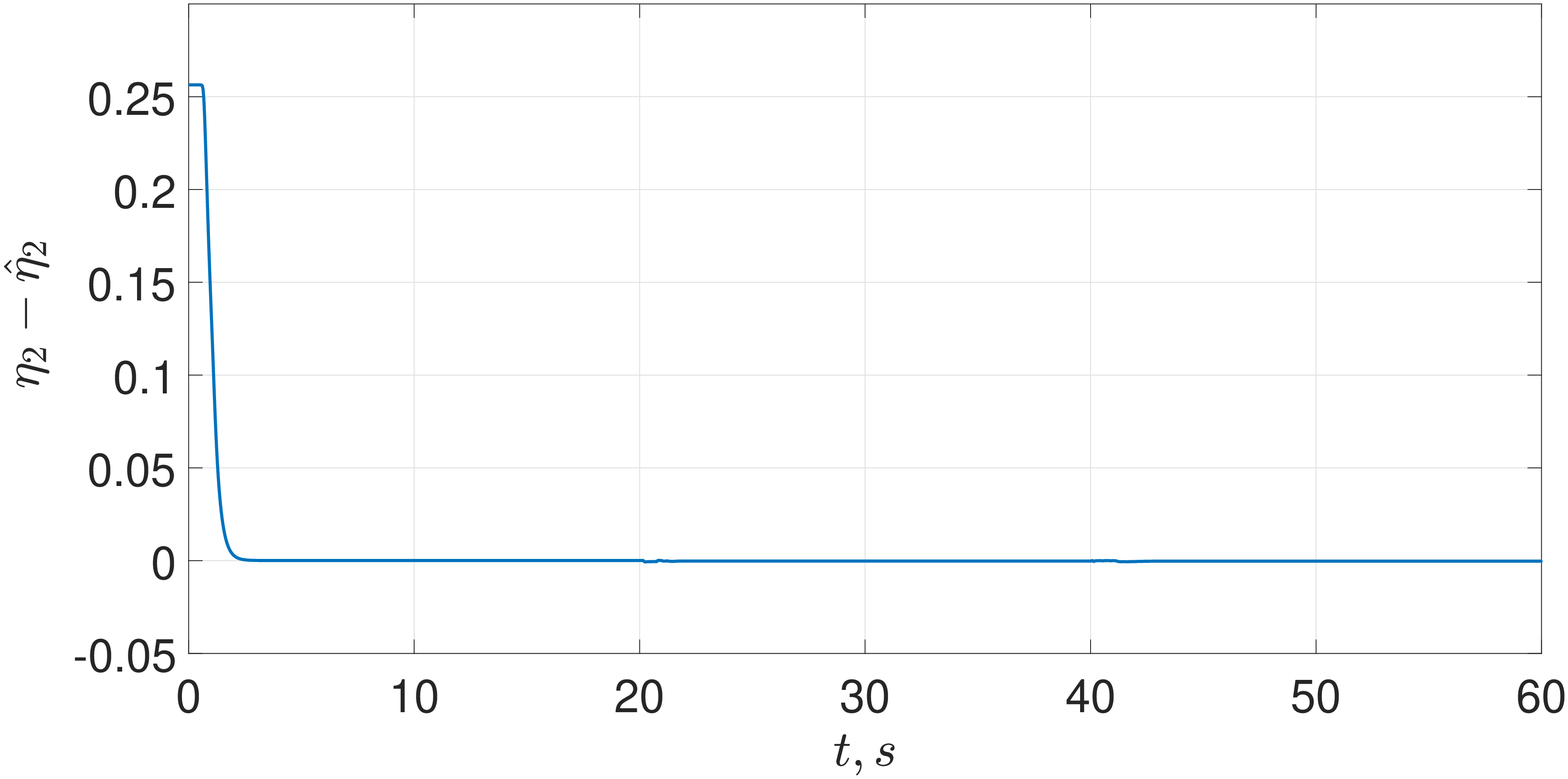}
	\caption{Transients of the error $\eta_2-\hat{\eta}_2$ ($\frac{b_2}{J}-\widehat{\frac{b_2}{J}}$)}
	\label{2_hat_eta_3}
\end{figure}

\begin{figure}[ht]
	\centering
	\includegraphics[width=1\linewidth]{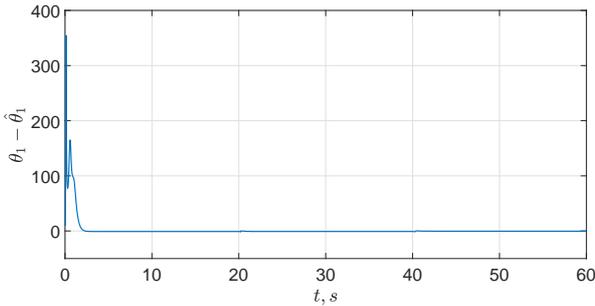}
	\caption{Transients of the error $\theta_1-\hat{\theta}_1$ ($\frac{k_1}{J}-\widehat{\frac{k_1}{J}}$)}
	\label{2_theta_1}
\end{figure}

If parameters $\theta_2$ and $\eta$ are known then we can estimate parameter $\theta_1$ using for instance standard gradietnt observer. We found paramete $\theta_1$ with standard gradient observer using $\theta_2$ and $\hat{\eta}$ instead real values with adaptation gain $10^6$ (see simulation result on the Fig. \ref{2_theta_1}).

Simulation results demonstrate convergence estimates to real values.

	\section{Concluding Remarks}
	\lab{sec6}
	%%%%%%%%%%%%%%%%%%
	%
	We have presented in this paper a constructive procedure to design GEC estimators for the parameters of two classes of nonlinear, NLP systems containing nonseparable nonlinearities of the form $e^{\theta_i h_i(u,y)}$. Although this class of nonlinearities seems to be very particular, as discussed in the Introduction, it appears in many practical applications, including the two thoroughly studied in the paper, and  is not amenable for the application of the existing parameter estimation techniques. The design procedure consists of the construction---from the non-separable NLP containing the exponential term---a new separable NLPRE, for which we can apply the LS+DREM estimator of \citep{ORTROMARA}. It is important to underscore that, to the best of our knowledge, only this estimator is capable of dealing with this kind of NLPREs. Moreover, the excitation requirement needed to ensure GEC is the very weak condition of IE defined in Assumption {\bf A4}.
	
	We would like to bring to the readers attention that techniques similar to the ones proposed here have been recently applied by the authors to solve two currently very relevant practical applications. Indeed, in \citep{BOBetalscl} we solve the problem of estimation of the parameters of the power coefficient of windmill generators in off-grid operation. The mathematical model of this system is of the form
	$$
	\dot y =-y^3(\theta_1 y -\theta_2)e^{-\theta_3 y},
	$$
	with $\theta \in \rea^3$ unknown parameters. Also, in \citep{BOBetalaut} we proposed a GEC parameter estimator for photo-voltaic arrays, whose dynamic model is of the form
	\begalis{
		\dot x  = - \theta_1 x -\theta_2 e^{b x}+\theta_3-\theta_4 u, 
		y  =  x - \theta_5 u,
	}
	with $\theta \in \rea^5$ unknown parameters, and the state $x(t) \in \rea$ {\em unmeasurable}. Notice that none of these applications fits into the class of systems considered in the paper.
	%%%%%%%%%%%%%%%%%
	%\section*{Acknowledgments}
\bibliography{Ortega_4}             % bib file to produce the bibliography
                                                     % with bibtex (preferred)
                                                   
%\begin{thebibliography}{xx}  % you can also add the bibliography by hand

%\bibitem[Able(1956)]{Abl:56}
%B.C. Able.
%\newblock Nucleic acid content of microscope.
%\newblock \emph{Nature}, 135:\penalty0 7--9, 1956.

%\bibitem[Able et~al.(1954)Able, Tagg, and Rush]{AbTaRu:54}
%B.C. Able, R.A. Tagg, and M.~Rush.
%\newblock Enzyme-catalyzed cellular transanimations.
%\newblock In A.F. Round, editor, \emph{Advances in Enzymology}, volume~2, pages
%  125--247. Academic Press, New York, 3rd edition, 1954.

%\bibitem[Keohane(1958)]{Keo:58}
%\newblock \emph{Power and Interdependence: World Politics in Transitions}.
%\newblock Little, Brown \& Co., Boston, 1958.

%\bibitem[Powers(1985)]{Pow:85}
%\newblock Is there a way out?
%\newblock \emph{Harpers}, pages 35--47, June 1985.

%\bibitem[Soukhanov(1992)]{Heritage:92}
%\newblock \emph{{The American Heritage. Dictionary of the American Language}}.
%\newblock Houghton Mifflin Company, 1992.

%\end{thebibliography}

%\appendix
%\section{A summary of Latin grammar}    % Each appendix must have a short title.
%\section{Some Latin vocabulary}              % Sections and subsections are supported  
                                                                         % in the appendices.
  \appendix
  %%%%%%%%%%%%%%
  \section{Proof of Lemmas}
  \lab{appa}
  %%%%%%%%%%%%%%%
  %
  \subsection*{Proof of Lemma \ref{lem1}}
  %%%%%%%%%%%%
  %
  We make the key observation that the function $h_4$ verifies
  $$
  \dot h_4 = \nabla^\top h_4 [f_1 + f_2 \calg(\eta)],
  $$
  where we have used \eqref{dotx1}.
  
  To simplify the notation in the sequel we define 
  \begin{footnotesize}
  	\begequ
  \lab{h5etal}
  h_5:=\nabla^\top h_4 f_1,\; h_6:=f_2^\top\nabla h_4,\;h_7:= {y_1-h_1 \over h_3},\;h_8:= -{h_2 \over h_3}.
  \endequ
  \end{footnotesize}
  We observe that using this notation we can write
  \begequ
  \lab{doth4}
  \dot h_4 = h_5 + h^\top_6 \calg(\eta).
  \endequ
  From the definition of $y_1$ in \eqref{y1} we get 
  \begali{
  	\nonumber
  	e^{h_4\theta_1}=h_7+h_8 \theta_2 \dif \theta_1 \dot h_4 e^{h_4\theta_1}=\dot h_7+\dot h_8 \theta_2}
  \begali{
  	\nonumber
  	\dif \theta_1 \dot h_4 e^{h_4\theta_1}=\dot h_7+\dot h_8 \theta_2}
  \begali{
  	\nonumber
  	\Longleftrightarrow \theta_1 \dot h_4 ( h_7 +h_8 \theta_2) =\dot h_7+\dot h_8 \theta_2}
  \begali{
  	\nonumber
  	\Longleftrightarrow \dot h_7 = \dot h_4 h_7 \theta_1 + \dot h_4 h_8 \theta_1 \theta_2 - \dot h_8 \theta_2}
  \begali{ 
  	\nonumber
  	\Longrightarrow \dot h_7 = [h_5 + h^\top_6 \calg(\eta)] h_7 \theta_1 + [h_5 + h^\top_6 \calg(\eta)] h_8 \theta_1 \theta_2}
  \begali{
  	\nonumber
  	- \dot h_8 \theta_2}
  \begali{
  	\nonumber
  	\Longleftrightarrow \dot h_7 =\theta_1h_5 h_7- \theta_2\dot h_8 + \theta_1\theta_2h_5 h_8 +\theta_1\calg^\top (\eta)h_7 h_6}
  \begali{
  	\nonumber
  	+ \theta_1\theta_2 \calg^\top(\eta)h_8  h_6, 
  }
  where  we used \eqref{doth4} to get the fourth implication. To obtain from the last identity a measurable NLPRE we apply the standard filtering technique \citep{IOASUNbook,TAObook}. Toward this end, we fix a constant $\lambda>0$ and define the stable filter ${\lambda \over p + \lambda}$, where $p:={d \over dt}$. Applying this filter to the last equation above, and recalling the definitions \eqref{h5etal}, we get appropriate vectors $Y_I$ and  $\phi_I$ for the NLPRE \eqref{nlpre1}  completing the proof.\footnote{Notice that the term $\phi_{I,2}$ can be computed without differentiation via the {\em proper} filtering ${\lambda p\over p + \lambda}\Big({h_2 \over h_3}\Big)$.}
  \subsection*{Proof of Lemma \ref{lem2}}
  %%%%%%%%%%%%
  %
  The Jacobian of $\calw(\theta)$ is given as
  $$
  \nabla \calw(\theta)=\begmat{1 & 0 &  0_{1 \times n_\eta}\\ 0 & 1 &  0_{1 \times n_\eta}\\  \theta_2 &  \theta_1 &  0_{1 \times n_\eta}\\ \calg(\eta) &  0_{p_\eta \times 1}  & \theta_1 \nabla \calg(\eta)\\ \theta_2 \calg(\eta)   & \theta_1 \calg(\eta) & \theta_1 \theta_2  \nabla \calg(\eta)}.
  $$
  The symmetric part of the matrix $T_\calw \nabla\calw(\theta)$ takes the form
  \begin{footnotesize}
  	\begali{
  	\nonumber
 &T_\calw  {\nabla}\calw(\theta)+ [\nabla \calw(\theta)]^\top T_\calw^\top \\
 	\nonumber
 &=\begmat{2 \alpha I_2 & \vdots & \begmat{ \sign(\theta_1)\calg^\top(\eta) T_\calg^\top \\  0_{1 \times n_\eta}}\\ \cdots & \cdot & \cdots \\ \begmat{ \sign(\theta_1) T_\calg \calg(\eta) & 0_{n_\eta \times 1}} &  \vdots & |\theta_1|\{T_\calg \nabla \calg(\eta) + [\nabla \calg(\eta)]^\top T_\calg^\top\}  }.
}
  \end{footnotesize}
  
  Let us introduce the notation
  \begali{
  	\nonumber
  	B  := \begmat{ \sign(\theta_1)\calg^\top(\eta) T_\calg^\top \\ \cdots \\ 0_{1 \times n_\eta}},
  }
  \begali{
  	\nonumber
  	C := |\theta_1|\{T_\calg \nabla \calg(\eta) + [\nabla \calg(\eta)]^\top T_\calg^\top\}.
  }
 
  A simple Schur complement calculation proves that the matrix $T_\calw  {\nabla}\calw(\theta)+ [\nabla \calw(\theta)]^\top T_\calw^\top$ is positive definite if and only if 
  \begequ
  \lab{c}
  C > {1 \over 2\alpha} B^\top B.
  \endequ
  On the other hand, from Assumption {\bf A3} we have that $C \geq |\theta_1| \rho_\calg I_{n_\eta}>0$. From which we conclude that  \eqref{c} holds for sufficiently large $\alpha$, concluding the proof. 
  \subsection*{Proof of Lemma \ref{lem3}}
  %%%%%%%%%%%%
  %
  To simplify the notation in the sequel we define 
  \begequ
  \lab{f3f4}
  f_3:= {f_1 \over h_3} ,\;f_4:= {1 \over h_3}f_2.
  \endequ
  We observe that using this notation and \eqref{sys2} we get the following chain of implications
  \begali{
  	\nonumber
  	{1 \over h_3} \ddot x  = e^{h_4\theta_1}  + f_3+f^\top_4 \calg(\eta) +  {u \over h_3}}
  \begali{
  	\nonumber
  	\dif - {\dot h_3 \over h^2_3}  \ddot x+  {1 \over h_3}{d^3 x \over dt^3} =  \theta_1 \dot h_4 \Big( e^{h_4\theta_1}\Big) + \dot f_3 + \dot f^\top _4 \calg(\eta)}
  \begali{
  	\nonumber
  	-  {\dot h_3 \over h^2_3}u +  {\dot u \over h_3}}
  \begali{
  	\nonumber
  	\Longleftrightarrow - {\dot h_3 \over h^2_3}  \ddot x+  {1 \over h_3}{d^3 x \over dt^3} =   \theta_1 \dot h_4 \Big( {1 \over h_3} \ddot x - f_3 - f^\top_4 \calg(\eta) - {u \over h_3} \Big)}
  \begali{
  	\nonumber
  	+ \dot f_3 + \dot f^\top _4 \calg(\eta) -  {\dot h_3 \over h^2_3}u +  {\dot u \over h_3}}
 \begali{
 	\nonumber
 	\hhh - {\dot h_3 \over h_3}  \ddot x+ {d^3 x \over dt^3} =   \theta_1   \dot h_4 \Big( \ddot x - h_3 f_3 - h_3 f^\top_4 \calg(\eta) -  u \Big)}
  \begali{
 	\nonumber
  	& + h_3\dot f_3 + h_3 \dot f^\top _4 \calg(\eta) -  {\dot h_3  \over h_3}u+  \dot u.
  }
  Moving to the left hand side the term which are {\em independent} of the unknown parameters we obtain the key identity
  \begali{
  	\lab{keyide}
  	\nonumber
  	&{d^3 x \over dt^3} - {\dot h_3 \over h_3}  \ddot x - h_3\dot f_3 +  {\dot h_3  \over h_3}u -  \dot u = \theta_1 \dot h_4 \Big( \ddot x -{1 \over h_3} f_3 -  u \Big) \\
  	&+\calg^\top(\eta)  h_3\dot f_4 - \theta_1  \calg^\top   {\dot h_4 \over h_3} f_4,}
    
  To obtain from \eqref{keyide} a measurable NLPRE we apply the standard filtering technique with the stable second order filter ${\lambda^2 \over (p + \lambda)^2}$. We observe that, due to the measurement of $\dot x$, the term
  $$
  {\lambda^2 \over (p + \lambda)^2}{d^3 x \over dt^3}={\lambda^2 p^2 \over (p + \lambda)^2}\dot x,
  $$
  is computable without differentiation. The same argument can be applied to all the other terms, except ${\dot h_3 \over h_3}  \ddot x$, $\dot h_4 \ddot x$ and $h_3\dot f_4$, which involve the unmeasurable signal $\ddot x$. To overcome this problem we invoke the Swapping Lemma \citep[Lemma 3.6.5]{SASBODbook}
  \begalis{
  	{\lambda \over p + \lambda} (\dot h_4 \ddot x) &= {\lambda \over p + \lambda} (h_4' \dot x \ddot x) =\hal {\lambda \over p + \lambda} \Big[h_4'p(\dot x^2)\Big]\\
  	& =  \hal\Big[ h_4' {\lambda p\over p + \lambda} (x^2)-{\lambda \over p + \lambda}\Big(h_4^{''} \dot x {p \over p + \lambda}(\dot x^2) \Big)\Big],
  }
  where the last right hand term can be computed without differentiation. Clearly, the same procedure can be applied to the term ${\dot h_3 \over h_3}  \ddot x$, leading to 
  \begalis{
  	&{\lambda \over p + \lambda} \Big({\dot h_3 \over h_3}  \ddot x \Big)  = {\lambda \over p + \lambda} ({h'_3 \over h_3} \dot x \ddot x) =\hal {\lambda \over p + \lambda} \Big[{h'_3 \over h_3}  p(\dot x^2)\Big]\\
  	& =  \hal\Big\{ {h'_3 \over h_3} {\lambda p\over p + \lambda} (x^2)+{\lambda \over p + \lambda}\Big[\Big( {(h_3')^2 \over h_3^2} - {h_3^{''} \over h_3}\Big) \dot x{p \over p + \lambda}(\dot x^2) \Big]\Big\}.
  }
  where, again,  the last right hand term can be computed without differentiation.
  
  Now, regarding the term $h_3\dot f_4$, from \eqref{f3f4}, we have that 
  \begalis{
  	h_3\dot f_4 & =\dot f_2 - {\dot h_3 \over h_3}f_2\\
  	& =\nabla_x f_2  \dot x +\nabla_{\dot x} f_2  \ddot x - {\dot h_3 \over h_3}f_2\\
  	& = \nabla_x f_2 \dot x +  \psi_a(x)\psi_b(\dot x) \ddot x  - {\dot h_3 \over h_3}f_2\\ 
  	& = \nabla_x f_2 \dot x +  \psi_a(x)\dot \psi_c(\dot x)  - {\dot h_3 \over h_3}f_2,
  }
  where we used Assumption {\bf A3} in the third identity and defined the function $\psi_c:=\int \psi_b(s)ds$. Applying the Swapping Lemma we can take care of the term $\psi_a\dot \psi_c$   
  \begalis{
  	{\lambda \over p + \lambda} [\psi_a  p (\psi_c)] &= \psi_a {\lambda p\over p + \lambda} (\psi_c)-{\lambda \over p + \lambda}\Big(\psi_a^{'} \dot x {p \over p + \lambda}(\psi_c) \Big),
  }
  which is clearly computable.
  
  Applying the second order filter to \eqref{keyide} and invoking the calculations above we get, after lengthy but straightforward calculations, get appropriate vectors $Y_{II}$ and  $\phi_{II}$ for the NLPRE \eqref{nlpre2} completing the proof.
  %%%%%%%%%%%%%%                                                                                                                                
                                                                         
\end{document}